\documentclass{amsart}
\usepackage{enumerate, amssymb}

\def\pdfsyncstart{}
\def\pdfsyncstop{}
\def\bdi{\pdfsyncstop\begin{diagram}}
\def\edi{\end{diagram}\pdfsyncstart}

\usepackage{diagrams}
\diagramstyle[scriptlabels,small]

\def\fm{{\mathfrak m}}
\def\fa{{\mathfrak a}}
\def\fp{{\mathfrak p}}

\def\lto{\longrightarrow}

\def\Ext{{\rm Ext}}
\def\Hom{{\rm Hom}}
\def\C{{\mathbb C}}

\def\Z{{\mathbb Z}}
\def\cB{{\mathcal B}}
\def\cP{{\mathcal P}}

\def\anotimes{{\hat\otimes}}

\DeclareMathOperator{\depth}{depth}

\DeclareMathOperator{\HH}{\rm HH}

\DeclareMathOperator{\tHH}{{\rm \widehat{HH}}}
\def\hto{{\hookrightarrow}}

\theoremstyle{definition}
\newtheorem{defn}{Definition}[section]

\theoremstyle{plain}
\newtheorem{prop}[defn]{Proposition}
\newtheorem{theorem}[defn]{Theorem}
\newtheorem{lem}[defn]{Lemma}
\newtheorem{cor}[defn]{Corollary}

\theoremstyle{remark}
\newtheorem{rem}[defn]{Remark}
\newtheorem{rems}[defn]{Remarks}

\newtheorem{sit}[defn]{}

\def\bpro{\begin{prop}}
\def\epro{\end{prop}}

\def\bthm{\begin{theorem}}
\def\ethm{\end{theorem}}

\def\bdfn{\begin{defn}}
\def\edfn{\end{defn}}

\def\brem{\begin{rem}}
\def\erem{\end{rem}}

\def\brems{\begin{rems}}
\def\erems{\end{rems}}

\def\bsit{\begin{sit}}
\def\esit{\end{sit}}

\def\blem{\begin{lem}}
\def\elem{\end{lem}}

\def\bcor{\begin{cor}}
\def\ecor{\end{cor}}

\def\bdia{\begin{diagram}}
\def\edia{\end{diagram}}

\def\ba{\begin{array}}
\def\ea{\end{array}}

\def\bnum{\begin{enumerate}}
\def\enum{\end{enumerate}}

\def\bprop{\begin{prop}}
\def\eprop{\end{prop}}

\def\be{\begin{equation}}
\def\ee{\end{equation}}

\def\bproof{\begin{proof}}
\def\eproof{\end{proof}}

\begin{document}
\title{Power Series Rings and Projectivity}

\author{Ragnar-Olaf Buchweitz}
\address{Dept.\ of Math., University of
Tor\-onto, Tor\-onto, Ont.\ M5S 2E4, Canada}
\email{ragnar@math.utoronto.ca}

\author{Hubert Flenner}
\address{Fakult\"at f\"ur Mathematik der Ruhr-Universit\"at,
Universit\"atsstr.\ 150, Geb.\ NA 2/72, 44780 Bochum, Germany}
\email{Hubert.Flenner@rub.de}

\thanks{The authors were partly supported by NSERC grant
3-642-114-80 and by the DFG Schwerpunkt ``Global Methods in Complex
Geometry".}


\begin{abstract}
We show that a formal power series ring $A[[X]]$ over a noetherian ring $A$ is not a projective module unless $A$ is artinian. However, if $(A,{\mathfrak m})$ is any local ring, then $A[[X]]$ behaves like a projective module in the sense that $\Ext^p_A(A[[X]], M)=0$ for all ${\mathfrak m}$-adically complete $A$-modules. The latter result is shown more generally for any flat $A$-module $B$ instead of $A[[X]]$. We apply the results to the (analytic) Hochschild cohomology over complete noetherian rings.
\end{abstract}

\maketitle

\section{Introduction}

By a classical result of Baer \cite{Baer}, the ring $\Z[[X]]$ of formal power series is not a projective $\Z$-module. In this paper we show  that more generally for any noetherian ring 
$A$ of dimension at least $1$ the ring $A[[X]]$ is not a projective $A$-module, see Theorem \ref{main1}. On the other hand we will deduce in Theorem \ref{main2} that any flat module $B$ over a local ring $(A,\fm)$ behaves like a projective module with respect to $\fm$-adically complete modules $M$ in the sense that $\Ext^p_A(B,M)=0$ for all $p\ge 1$.

Our motivation for these results comes from Hochschild cohomology. If $A\to B$ is a flat morphism of (commutative) rings then the usual bar resolution $\cB_\bullet$ provides a flat resolution of $B$ as a module over the enveloping algebra $B^{(e)}:=B\otimes_AB$, and by definition Hochschild cohomology $\HH^p(B/A,M)$ is the cohomology of the complex $\Hom_{B^{(e)}}(\cB_\bullet, M)$. If $B$ is projective as an $A$-module then $\cB_\bullet$ is even a complex of projective $B^{(e)}$-modules and so  
$$
\HH^p(B/A,M)\cong \Ext^p_{B^{(e)}}(B, M)\,.
$$
In Proposition \ref{hoch}, we apply Theorem \ref{main2} to show that these isomorphisms still prevail, as long as $B$ is flat over $A$ and $M$ is a complete module over $B^{(e)}$.

Next assume that $A\to B$ is a flat morphism of complete local rings with isomorphic residue fields. Using complete tensor products one can introduce the completed Bar-resolution $\hat\cB_\bullet$ to define an analytic version of Hochschild cohomology $\tHH^p(B/A,M)$ as the cohomology of the complex $\Hom_{\hat B^{(e)}}(\hat \cB_\bullet, M)$. As an application of the preceding results we will show in Section 3 that in analogy with the isomorphisms above there are as well isomorphisms 
$$
\tHH^p(B/A,M)\cong \Ext^p_{\hat B^{(e)}}(B, M)
$$
for complete $\hat B^{(e)}$-modules $M$.

\section{Complete modules and projectivity}

By a result of Baer \cite{Baer} the ring $\Z[[X]]$ is not a projective $\Z$-module. Using a variant of his idea of proof we can show the following result.   

\bthm\label{main1}
If $A$ is a commutative noetherian ring of dimension $\ge 1$, then
$A[[X]]$ is not a projective $A$-module.
\ethm

\bproof
As $A$ is noetherian by assumption, every ideal $\fa$ of $A$ is finitely generated and so $A[[X]]/\fa A[[X]]\cong (A/\fa)[[X]]$. Moreover, if $A[[X]]$ is projective as $A$-module so is $A[[X]]/\fa A[[X]]$ as a module over $A/\fa$. Replacing $A$ by $A/\fp$ for a minimal prime $\fp$ of $A$ we may, and will, henceforth assume that $A$ is a domain.  

Now fix a maximal ideal $\fm$ of $A$ and consider the subset
$Z\subseteq A[[X]]$ of all power series $F=\sum a_iX^i$ with
$\lim_{i} a_i=0$ in the $\fm$-adic topology. Clearly, $Z$ is an
$A$-submodule of $A[[X]]$ that contains the polynomial ring
$A[X]\subseteq A[[X]]$. We first show that
$$
Z=A[X]+\fm Z.\leqno (*)
$$
Indeed, if $F=\sum a_iX^i$ is in $Z$ then there is a sequence
of numbers $n_i$ with $a_i\in\fm^{n_i}$ and $\lim n_i=\infty$.
We may assume that $n_i\ge 1$ for $i\ge k$.  Thus, if
$f_1,\ldots,f_n$ is a system of generators of
$\fm$ as an $A$-module, then for $i\ge k$ we can write
$a_i=\sum_{j=1}^n a_{ij}f_j$ with $a_{ij}\in\fm^{n_i-1}$, whence 
$$
F=\sum_{i=0}^{k-1}a_iX^i+\sum_{j=1}^n F_jf_j\,,
$$
where $F_j:=\sum_{i\ge k}a_{ij}X^i$ is again in $Z$.

Now assume $A[[X]]$ were a projective $A$-module. It is then in particular contained in a free $A$-module $A^{(I)}$. Restricting the inclusion to $Z$ will then yield
$$
Z\hto A^{(I)}.
$$
As $A[X]$ is free on a countable basis, the image of $A[X]$ is
contained in some direct summand
$A^{(J)}$ of $A^{(I)}$, where $J\subseteq I$ is a countable subset. With $M:=Z/A[X]$,
there is hence an induced $A$--linear map
$$
\varphi: M \lto A^{(I\backslash J)}\,.
$$
By $(*)$ above, $\fm M=M$, whereas $\bigcap_{k\ge 0}\fm^kA^{(I\backslash
J)}=0$, whence $\varphi$ is necessarily the zero map and so already $Z\subseteq A^{(J)}$.

On the other hand, if $t\in\fm$ is a nonzero element, then the $A$-linear map
$$
A[[X]]\to Z\quad\mbox{with}\quad F=\sum a_iX^i\mapsto
\sum a_it^iX^i
$$
is injective, whence $A[[X]]$ can as well be realized as a submodule of $A^{(J)}$. We now show that this is impossible.

Indeed, if $A$ is uncountable, then the formula for Vandermonde's determinant yields that the power series 
$$
F_a:=\sum_{i=0}^\infty a^iX^i\,,\quad a\in A\setminus\{0\}\,,
$$
form $A$-linearly independent elements of $A[[X]]$. Hence, with $K$ the field of fractions of $A$, the vector space $A[[X]]\otimes_AK$ has uncountable dimension over $K$ contradicting the fact that it can be embedded into the countably generated vector space $A^{(J)}\otimes_AK\cong K^{(J)}$. Finally, if $A$ is countable then so is $A^{(J)}$, while $A[[X]]$ is uncountable using Cantor's argument. Hence we obtain the desired contradiction concluding the proof that $A[[X]]$ cannot be projective.
\eproof

\brems
1. A minor variant of the argument shows as well that for an analytic algebra $A$ over a valued field, the convergent power series ring $A\{X\}$ is neither projective as an $A$-module. For $A=\C\{t\}$, this was already observed by Wolffhardt
\cite[Satz 9]{Wo}. 

2. As pointed out by Avramov, if $R$ is a non-local domain and if $\hat R$ is its
completion with respect to some maximal ideal $\fm\subseteq R$,
then $\hat R$ is as well not projective as an $R$-module. In fact, if
$a\in R\backslash \fm$ is not a unit then $R/aR\ne 0$ whereas $\hat
R\otimes_R(R/aR)\cong \hat R/a\hat R=0$. Hence $\hat R$ is not faithful as an $R$-module and so cannot be projective. 

This shows in particular directly that $A[[X]]$, or $A\{X\}$ in the analytic case, is never a projective $R=A[X]$--module.
\erems

As was shown in \cite{Baer}, the group $\Ext^1_\Z(\Z[[X]], T)$ vanishes whenever $T$ is a finitely generated torsion $\Z$-module or, even more generally, $T$ is of bounded torsion. We will supplement this observation by the following fact.

\bthm\label{main2}
Let $A$ be a ring and $\fm\subseteq A$ a maximal ideal. If $B$ is a flat $A$-module then 
$\Ext^p_A(B, M)= 0$ for all $p\ge 1$ and each $\fm$-adically complete $A$-module $M$.
\ethm

We first recall the following well-known result; see \cite[II.3.Cor.2 to Prop.5]{AC}; and give its simple proof.
\blem\label{lem1}
Every flat module is free over a local ring $(A,\fm)$ whose maximal ideal is nilpotent, that is, $\fm^n=0$ for some $n\ge 1$. 
\elem

\bproof
The lemma is certainly true if $\fm=0$, as then $A=A/\fm$ is a field. In the general case, consider a flat $A$-module $B$ and elements $F_i\in B$, $i\in I$, whose residue classes form a basis of $B/\fm B$ as a vector space over $A/\fm$. Let us show that these elements then form a basis of $B$ as an $A$-module. If $B'$ is the submodule generated by these elements then $B/B'=\fm \cdot (B/B')$ and so 
$$
B/B'=\fm\cdot B/B'=\fm^2\cdot B/B'=\cdots =\fm^n\cdot B/B'=0.
$$
Thus, the elements $F_i$ generate $B$, and we now show that they are also linearly independent. For $B$ is flat, and, if $K=\ker(A^{(I)}\to B)$ denotes the kernel of the map defined by the $F_i$, tensoring the exact sequence 
$$
0\to K \to A^{(I)}\to B\to 0
$$
with $A/\fm$ over $A$ results in the exact sequence 
$$
0\to K/\fm K \to (A/\fm)^{(I)}\xrightarrow{\ \cong\ } B/\fm B\to 0\,.
$$
This shows first $K/\fm K=0$ and then the same argument as before yields $K=0$. Thus $B\cong A^{(I)}$ is free over $A$, as claimed.
\eproof

\blem\label{lem2}
Let $A$ be a ring and $\fm\subseteq A$ a maximal ideal. If $B$ is a flat $A$-module, then 
$\Ext^p_A(B, M)= 0$ for all $p\ge 1$ and every $A$-module $M$ with $\fm^n M=0$ for some $n\ge 0$.
\elem

\bproof
As $B$ is flat over $A$, we have 
$$
\Ext^p_A(B,M)\cong \Ext^p_{A/\fm^n }(B/\fm^n B, M)\quad \forall p\ge 0\,. \leqno(*)
$$
Recall the simple argument: if $F_\bullet$ is a projective resolution of $B$ as an $A$-module then by flatness of $B$ the complex $F_\bullet/\fm^n  F_\bullet$ is a projective resolution of $B/\fm^n B$ as an $A/\fm^n $-module. Hence $\Ext^p_{A/\fm^n }(B/\fm^n  B, M)$ can be computed from the complex 
$$
\Hom_{A/\fm^n } (F_\bullet/\fm^n F_\bullet, M)\cong \Hom_A(F_\bullet, M)\,.
$$
As the complex on the right has cohomology $\Ext^p_A(B,M)$, the claim $(*)$ follows. 
Flatness of $B$ over $A$ implies that $B/\fm^n B$ is flat as an $A/\fm^n $-module, whence the lemma is  a consequence of $(*)$ and Lemma \ref{lem1}. 
\eproof

Finally, recall the following simple fact. 

\blem\label{lem3}
Let $A$ be a ring and $\{f_n:H_{n+1}\to H_n\}_{n\ge 0}$ an inverse system of $A$-modules with surjective  transition maps $f_n$. Setting $H:=\prod_{n\ge 0}H_n$, the map $f:H\to H$ with $(h_n)_{n}\mapsto (h_n-f_n(h_{n+1}))_{n}$ is surjective and has kernel $\displaystyle\lim_\leftarrow H_n$.  \qed
\elem

Now we turn to the
\medskip

\noindent
{\it Proof of Theorem {\ref{main2}}.}
For $M_n:=M/\fm^{n+1}M$, Lemmata  \ref{lem1} and \ref{lem2} give that
$$
\Ext^p_A(B,M_n)=0\quad\text{for all}\quad p\ge 1\,.
$$
Taking the direct product $H:=\prod_{n\ge 0}M_n$ this implies
$$
\Ext^p_A(B, H)=0\quad\mbox{ for all }p\ge 1\,,
$$
since the formation of $\Ext$ is compatible with direct products in the second component. As $M$ is complete as an $A$-module there is an exact sequence 
$$
0\to M\to  H=\prod_{n\ge 0}M_n \xrightarrow{f} H=\prod_{n\ge 0}M_n \to 0\,,
$$
see Lemma \ref{lem3}.
Applying $\Hom_A(B,-)$ gives the long exact $\Ext$-sequence
$$
0\to \Hom_A(B,M) \to \Hom(B, H) \xrightarrow{f^{*}} \Hom_A(B,H) \to \Ext^1_A(B,M)\to\cdots
$$
As $\Ext^p_A(B, H)=0$ for $p\ge 1$ it follows that $\Ext^p_A(B,M)=0$ for $p\ge 2$. It remains to show that $\Ext^1_A(B, M)=0$ or, equivalently, that the map
$$
f^{*}: \Hom_{A}(B, H)\cong \prod_{n\ge 0}\Hom_A(B, M_n) \to \Hom_A(B,H)\cong \prod_{n\ge 0}\Hom_A(B, M_n) 
$$ 
is surjective. However, this is immediate from Lemma \ref{lem3} as the transition maps
$$
\Hom_A(B, M_{n+1}) \to \Hom_A(B, M_n) 
$$
are all surjective by Lemma \ref{lem2}.

\bcor\label{cor}
If $A$ is a complete local noetherian ring and $B$ is a flat $A$-module then $\Ext^p_A(B,M)=0$ for all $p\ge 1$ and all finite $A$-modules $M$.\qed
\ecor

\begin{rems}
1. The corollary itself is not new, indeed Jensen established in \cite[Thm.8.1]{Je2}: 
If $A$ is commutative noetherian then it is a product of a finite number of complete 
local rings if and only if $\Ext_A^i(B,M)=0$, for $i\geq 1$, whenever $B$ is flat and 
$M$ is finite over $A$.

2.  In a similar vein, Frankild establishes in \cite[Cor.3.7]{Fr}: For a local noetherian 
ring $(A,\fm)$, an $A$--module $B$ of finite projective dimension 
and an $A$--module $M$ that is $\fm$-adically complete, $\Ext_{A}^{i}(B,M) = 0$ 
for $i > \depth A -\depth B$. 

As any flat $A$--module is of finite projective dimension; see \cite[p.164]{Je1}; 
with its depth either equal to the depth of the ring or infinite, this result specializes 
to Corollary \ref{cor}.

3.
For non-complete modules $M$, the problem as to whether $\Ext^p_A(B,M)$ 
vanishes for $p\ge 0$ is much more intricate. For instance, the famous Whitehead 
problem asked whether the vanishing of $\Ext^1_\Z(B,\Z)$ implies that the module 
$B$ is free. As was shown by Shelah, this depends on the model of set theory used, 
see \cite{Sh}. 

4.
We do not know whether a result similar to Corollary \ref{cor} also holds in the 
analytic category. More precisely, let $A\to A'$ be a homomorphism of analytic 
algebras over a valued field and $B$ a finite $A'$-module that is flat over $A$. 
It is natural to ask whether then
$\Ext^p_A(B,M)=0$ for all finite $A$-modules $M$ and all $p\ge 1$.
\end{rems}

\section{Hochschild cohomology of complete algebras}

Let us first recall the definition of Hochschild cohomology of a morphism of 
commutative rings $A\to B$. Denoting $B^{\otimes n}$ the usual $n$-fold 
tensor product over $A$, the bar resolution 
$$
\cB_\bullet:\quad \cdots B^{\otimes n} \to \cdots\to  B^{\otimes 2}\to B\to 0
$$
provides a resolution of $B$ as a $B^{(e)}:=B^{\otimes 2}$-module, see \cite{CE} 
or \cite{Lo}. Note that this resolution is flat, respectively projective over $B^{(e)}$ 
provided $B$ has the same property over $A$. For any  $B^{(e)}$-module $M$ 
the modules
$$
\HH_p(B/A,M)=H_p(\cB_\bullet\otimes_{B^{(e)}}M)
$$
and 
$$
\HH^p(B/A,M)=H^p(\Hom_{B^{(e)}}(\cB_\bullet,M))
$$
are called, repectively, the Hochschild homology and cohomology of $B/A$ 
with values in $M$. 

If $\cP_{\bullet}\to B\to 0$ is a projective resolution of $B$ over $B^{(e)}$, 
then there exists a comparison map $\cP_{\bullet}\to \cB_{\bullet}$ over the 
identity of $B$ that is a homomorphism of complexes and unique up to homotopy. Accordingly, for any $B^{(e)}$-module $M$ there are homomorphisms
$$
\HH^p(B/A,M)\to \Ext^p_{B^{(e)}}(B, M)\,,
$$
that are functorial in $M$. These maps are isomorphisms as soon as 
$B$ is a projective $A$--module, but an application of Theorem \ref{main2} 
shows that these maps are also isomorphisms in the following wider context.
\bpro
\label{hoch}
If $A\to B$ is a flat homomorphism of commutative rings and $M$ a 
$B^{(e)}$-module that is $\fm$-adically complete for some maximal ideal 
$\fm\subseteq B^{e}$, then the natural maps 
$\HH^p(B/A,M)\to \Ext^p_{B^{(e)}}(B, M)$ are isomorphisms for each $p\ge 0$.
\epro

\begin{proof}
Consider the pair of spectral sequences with the same limit
\begin{align*}
{}'E_{1}^{pq}= \Ext^{q}_{B^{(e)}}(\cB_{p}, M)&\Longrightarrow
{\mathbb E}^{p+q}\quad\text{and}\\
{}''E_{2}^{pq} = \Ext^{p}_{B^{(e)}}(H_{-q}(\cB_{\bullet}), M)&\Longrightarrow
{\mathbb E}^{p+q}\,.
\end{align*}
As $\cB_{\bullet}$ is a resolution of $B$, the second of these spectral sequences 
degenerates and identifies the limit as ${\mathbb E}^{n}\cong \Ext^{n}_{B^{(e)}}(B, M)$. Concerning the first spectral sequence, each term $\cB_{p}$ is a flat 
$B^{(e)}$-module, and so Theorem \ref{main2} shows that ${}'E_{1}^{pq}= 0$ 
for $q\neq 0$ and any $p$. Hence the first spectral sequence degenerates as 
well, identifying the limit as ${\mathbb E}^{n}\cong \HH^{n}(B/A, M)$.
\end{proof}
For the rest of this section, let us assume that $A\to B$ is a homomorphism of 
complete local rings such that the extension of residue fields $A/\fm_A\to B/\fm_B$ 
is an isomorphism. In the usual construction of Hochschild (co-)homology, 
as sketched above, one can then replace the ordinary tensor product by the 
complete one, $B^{\anotimes n}$, that is the completion of $B^{\otimes n}$ 
with respect to the maximal ideal 
$\ker( B^{\otimes n}\to (B/\fm_B)^{\otimes n}\cong B/\fm_B)$. The result is 
the analytic Bar-resolution $\hat\cB_\bullet$ that gives rise to analytic 
Hochschild (co-)homology
$$
\tHH_p(B/A,M)=H_p(\hat\cB_\bullet\otimes_{\hat B^{(e)}}M)
$$
and 
$$
\tHH^p(B/A,M)=H^p(\Hom_{\hat B^{(e)}}(\hat\cB_\bullet,M)). 
$$
There is obviously a canonical morphism $\cB_{\bullet}\to \hat \cB_{\bullet}$, 
whence there are induced maps
$$
\alpha_p:\HH_p(B/A,M)\to \tHH_p(B/A,M)
$$
and 
$$
\beta^p: \tHH^p(B/A,M)\to \HH^p(B/A,M). 
$$
As an application of the results of section 2 we obtain the following proposition that was actually the motivation for this paper.

\bpro\label{hochschild}
If $B$ is a flat $A$-algebra and $M$ is a complete ${\hat B^{(e)}}$-module, then there is a natural isomorphism
$$
 \tHH^p(B/A,M)\cong \Ext^p_{\hat B^{(e)}}(B,M)
$$
for each $p\ge 0$. 
\epro

\bproof
The analytic Bar-resolution $\hat\cB_\bullet$ provides a flat resolution of $B$ as a $\hat B^{(e)}$-module. Using again Theorem \ref{main2} and the same spectral sequence argument as above shows that the cohomology groups of the complex
$$
\Hom_{\hat B^{(e)}}(\hat\cB_\bullet,M)
$$
can be identified with the $\Ext$-groups $\Ext^p_{\hat B^{(e)}}(B,M)$ as claimed. 
\eproof

\brems
1. As is well known, the map $\alpha_p$ is not an isomorphism, in general. 
For instance, $\HH_1(B/A, B)\cong D_A(B)=\Omega^{1}_{B/A}$ is the 
{\em universal\/} module of differentials for $A\to B$ that is not finitely 
generated, in general. On the other hand, $\tHH_1(B/A, B)\cong D^f_A(B)
={\widehat\Omega}^{1}_{B/A}$ is the {\em universally finite\/} module of 
differentials, which is always a finite module (see, e.g., \cite{SS} or \cite{EGA IV}). 

2. We do not know whether the maps $\beta^p$ are isomorphisms when 
$A\to B$ is flat and $M$ is complete. For this it would be necessary to know 
whether $\cB_\bullet\otimes_{B^{(e)}}\hat B^{(e)}$ is a resolution of $B$. 
However, $B^{(e)}$ is not noetherian, in general, whence it is not clear 
whether $B^{(e)}\to\hat B^{(e)}$ is always a flat map. 
\erems

{\em Acknowledgement:\/}
We thank L.\ Avramov and S.\ Iyengar for valuable comments and pointers to the literature.

\end{document}